\documentclass[amstex,12pt,english,amssymb]{article}

\usepackage{mathtext}
\usepackage[cp1251]{inputenc}
\usepackage[T2A]{fontenc}
\usepackage[english]{babel}
\usepackage[dvips]{graphicx}
\usepackage{amsmath}
\usepackage{amssymb}
\usepackage{amsxtra}
\usepackage{latexsym}
\usepackage{ifthen}

\begin{document}

\newcounter{lemma}[section]
\newcommand{\lemma}{\par \refstepcounter{lemma}%
{\bf Lemma \arabic{section}.\arabic{lemma}.}}
\renewcommand{\thelemma}{\thesection.\arabic{lemma}}

\newcounter{corollary}[section]
\newcommand{\corollary}{\par \refstepcounter{corollary}%
{\bf Corollary \arabic{section}.\arabic{corollary}.}}
\renewcommand{\thecorollary}{\thesection.\arabic{corollary}}

\newcounter{remark}[section]
\newcommand{\remark}{\par \refstepcounter{remark}%
{\bf Remark \arabic{section}.\arabic{remark}.}}
\renewcommand{\theremark}{\thesection.\arabic{remark}}

\newcounter{theorem}[section]
\newcommand{\theorem}{\par \refstepcounter{theorem}%
{\bf Theorem \arabic{section}.\arabic{theorem}.}}
\renewcommand{\thetheorem}{\thesection.\arabic{theorem}}

\newcounter{proposition}[section]
\newcommand{\proposition}{\par \refstepcounter{proposition}%
{\bf Proposition \arabic{section}.\arabic{proposition}.}}
\renewcommand{\theproposition}{\thesection.\arabic{proposition}}

\newcommand{\proof}{{\it Proof.\,\,}}

\renewcommand{\theequation}{\arabic{section}.\arabic{equation}}

\def\Xint#1{\mathchoice
   {\XXint\displaystyle\textstyle{#1}}%
   {\XXint\textstyle\scriptstyle{#1}}%
   {\XXint\scriptstyle\scriptscriptstyle{#1}}%
   {\XXint\scriptscriptstyle\scriptscriptstyle{#1}}%
   \!\int}
\def\XXint#1#2#3{{\setbox0=\hbox{$#1{#2#3}{\int}$}
     \vcenter{\hbox{$#2#3$}}\kern-.5\wd0}}
\def\dashint{\Xint-}
\author{Ryazanov V.  and  Sevost'yanov E.}

\title{ %{\leftline{\small UDC 517.5}}
\medskip Equicontinuity of mappings \\ quasiconformal in the mean}
\maketitle
\begin{abstract}
It is stated equicontinuity and normality of families
$\frak{R}^{\Phi}$ of the so--called ring $Q(x)$--homeomorphisms with
integral constraints of the type
$\int\Phi\left(Q(x)\right)dm(x)<\infty$ in a domain $D\subset{\Bbb
R}^n,$ $n\ge 2.$ It is shown that the found conditions on the
function $\Phi$ are not only sufficient but also necessary for
equicontinuity and normality of such families of mappings. It is
also given applications of these results to families of mappings in
the Sobolev class $W_{loc}^{1,n}.$
\end{abstract}

\section{Introduction}
Here $dm(x)$ corresponds to the Lebesgue measure in a domain $D$ in
${\Bbb R}^n,$ $n\ge 2$.\medskip

In the theory of  mappings called quasiconformal in the mean,
conditions of the type
\begin{equation}\label{eq2} \int\limits_{D} \Phi
(Q(x))\ dm(x)\  <\ \infty\end{equation} are standard for various
characteristics $Q$ of these mappings, see e.g. \cite{Ah},
\cite{Bi}, \cite{Gol}, \cite{GMSV}, \cite{Kr$_1$}--\cite{Ku},
\cite{Per}, \cite{Pes}, \cite{Rya} and \cite{UV}. The study of
classes with the integral conditions (\ref{eq2}) is also actual in
the connection with the recent development of the theory of
degenerate Beltrami equations and the so--called mappings with
finite distortion, see e.g. related references  in the monographs
\cite{IM$_1$} and \cite{MRSY}.\medskip

In the present paper we study the problems of equicontinuity and
normality for wide classes of the so--called ring
$Q(x)-$homeomorphisms with the condition (\ref{eq2}) and give the
corresponding applications to Sobolev's classes.\medskip

Recall that the {\bf (conformal) modulus} of a family $\Gamma$ of
curves $\gamma$ in ${\Bbb R}^n$, $n\ge 2$, is the quantity
$$
M(\Gamma)=\inf_{\rho \in \,{\rm adm}\,\Gamma} \int\limits_{{\Bbb
R}^n} \rho ^n (x)\ \ dm(x)
$$
where a Borel function $\rho:{\Bbb R}^n\,\rightarrow [0,\infty]$ is
{\bf admissible} for $\Gamma$, write $\rho \in {\rm adm} \,\Gamma $,
if
$$
\int\limits_{\gamma}\rho (x)\ \ |dx|\ge 1\ \ \ \ \ \ \ \forall\
\gamma \in \Gamma\ .
$$

One of the equivalent geometric definitions of {\bf
$K-$quasiconformal mappings} $f$ with $K\in [1,\infty)$ given in a
domain $D$ in ${\Bbb R}^n,$ $n\ge 2,$ is reduced to the inequality
\begin{equation}\label{eq4*}
M(f\Gamma)\le K\,M(\Gamma)
\end{equation}
that holds for an arbitrary family  $\Gamma$  of curves $\gamma$ in
the domain $D$. \medskip

Similarly,  given a domain $D$ in ${\Bbb R}^n,$ $n\ge 2,$  and a
(Lebesgue) measurable function $Q:D\to[1,\infty]$, a homeomorphism
$f:D\to{\overline{{\Bbb R}^n}}$, ${\overline{{\Bbb R}^n}}={\Bbb
R}^n\cup\{\infty\}$, is called {\bf $Q(x)$ -- homeomorphism} if
\begin{equation} \label{eq19*}
M(f\Gamma )\le \int\limits_D Q(x)\cdot \rho^n (x)\ \ dm(x)
\end{equation}
for every family  $\Gamma$  of curves $\gamma$ in $D$ and every
$\rho \in {\rm adm} \,\Gamma $, see e.g. \cite{MRSY}.\medskip

In the case $Q(x)\le K$ a.e., we again come to the inequality
(\ref{eq4*}). In the general case, the latter inequality means that
the conformal modulus of the family $f\Gamma$ is estimated by the
modulus  $M_Q$ of $\Gamma$ with the weight $Q,$
$M(f\Gamma)\le M_{Q}(\Gamma),$
see e.g. \cite{AC$_1$}. The inequality of the type (\ref{eq19*}) was
first stated by O. Lehto and K. Virtanen for quasiconformal mappings
in the plane, see Section V.6.3 in \cite{LV}. The relation of the
type (\ref{eq19*}) was also stated by K. Bishop, V. Gutlyanskii, O.
Martio and M. Vuorinen in \cite{BGMV} for quasiconformal mappings in
space where $Q(x)$ is equal to $K_I(x,f).$\medskip

Recall that the {\bf inner dilatation} of a mapping $f:D\rightarrow
{\Bbb R}^n,$ $n\ge 2,$ at a point $x\in D$ of differentiability for
$f$ is
$$K_I(x,f)\,\,=\,\,\frac{|J(x,f)|}{{l\left(f^{\,\prime}(x)\right)}^n}$$
%
%\end{equation}
%
if $J(x,f) \ne 0,$\,\, $K_I(x,f)=1$ if  $f^{\,\prime}(x)=0,$ and
$K_I(x,f)=\infty $ at the rest points, where $J(x,f)$ is the
Jacobian of $f$ at $x$ and
$$l\left(f^{\,\prime}(x)\right)\,=\,\,\,\inf\limits_{h\in {\Bbb R}^n
\backslash \{0\}} \frac {|f^{\,\prime}(x)h|}{|h|}\,\,.$$\medskip

The following notion generalizes and localizes the above notion of a
$Q$--ho\-me\-o\-mor\-phism. It is motivated by the ring definition
of Gehring for qua\-si\-con\-for\-mal mappings, see e.g.
\cite{Ge$_3$}, introduced first in the plane, see \cite{RSY$_3$},
and extended  later on to the space case in \cite{RS}, see also
Chapters 7 and 11 in \cite{MRSY}. Let $E,$ $F\subset\overline{{\Bbb
R}^n}$ be arbitrary sets. Denote by $\Gamma(E,F,D)$  a family of all
curves $\gamma:[a,b]\rightarrow \overline{{\Bbb R}^n}$ joining $E$
and $F$ in $D,$ i.e. $\gamma(a)\in E,\gamma(b) \in F$ and
$\gamma(t)\in D$ as $t \in (a, b).$\medskip

Given a domain $D$ in ${\Bbb R }^n,$ $n\ge 2,$ a (Lebesgue)
measurable function $Q:D\rightarrow\,[0,\infty]$, $x_0\in D,$
a homeomorphism $f:D\rightarrow \overline{{\Bbb R}^n}$ is said to be
a {\bf ring $Q$--homeomorphism at the point $x_0$} if
\begin{equation}\label{eq1}
M\left(f\left(\Gamma\left(S_1,\,S_2,\,R\right)\right)\right)\ \le
\int\limits_{R} Q(x)\cdot \eta^n(|x-x_0|)\ dm(x)
\end{equation}
for every ring $R=R(r_1,r_2,x_0)$ $=\{ x\,\in\,{\Bbb R}^n :
r_1<|x-x_0|<r_2\} $ and the spheres $S_i=S(x_0, r_i)=\{x\in {\Bbb
R}^n: |x-x_0|=r_i\}$, where $0<r_1<r_2< r_0\,\colon =\,{\rm dist}\,
(x_0,\partial D),$ and every measurable function $\eta :
(r_1,r_2)\rightarrow [0,\infty ]$ such that
$$\int\limits_{r_1}^{r_2}\eta(r)\ dr\ \ge\ 1\,.$$
$f$ is called a {\bf ring $Q$--homeomorphism in the domain} $D$ if
$f$ is a ring $Q$--ho\-meo\-mor\-phism at every point $x_0\in D$.
Note that, in particular, homeomorphisms $f:D\rightarrow
\overline{{\Bbb R}^n}$ in the class $W_{loc}^{1,n}$ with
$K_I(x,f)\in L_{loc}^1$ are ring $Q$--ho\-me\-o\-mor\-phisms with
$Q(x)=K_I(x,f),$ see e.g. Theorem 4.1 in \cite{MRSY}.
\medskip

The notion of ring $Q$--ho\-me\-o\-mor\-phism can be extended in the
natural way to $\infty$.  More precisely, under $\infty\in D
\subseteq \overline{{\Bbb R}^n}$ a homeomorphism $f:D\rightarrow
\overline{{\Bbb R}^n}$  is called a {\bf ring
$Q$--ho\-me\-o\-mor\-phism at ${\bf \infty}$} if the mapping
$\widetilde{f}=f\left(\frac{x}{\,|x|^2}\right)$  is a ring
$Q^{\,\prime}$--ho\-me\-o\-mor\-phism at the origin with
$Q^{\,\prime}(x)=Q\left(\frac{x}{\,|x|^2}\right).$ In other words, a
mapping $f:{\Bbb R}^n\rightarrow \overline{{\Bbb R}^n}$ is a ring
$Q$--ho\-me\-o\-mor\-phism at $\infty$ iff
$$M\left(f\left(\Gamma\left(S(R_1), S(R_2), R\right)\right)\right)\le
\int\limits_{R} Q(y)\cdot \eta^n\left(|y|\right)dm(y)$$
%\end{equation}
holds for every ring $R=R(R_1, R_2, 0)=\{y\in {\Bbb R}^n:
R_1<|y|<R_2\}$ in $D$ with $0<R_1<R_2<\infty,$ $S(R_i)=\{x\in {\Bbb
R}^n: |x|=R_i\}$ and for every measurable function $\eta :
(R_1,R_2)\rightarrow [0,\infty ]$ with
%
%
%\begin{equation}\label{eq11}
%
$\int\limits_{R_1}^{R_2}\eta(r)\ dr\ \ge\ 1\,.$

\section{Preliminaries}
\setcounter{equation}{0}

Let $(X,d)$ and $\left(X^{\,\prime}, d^{\,\prime}\right)$ be metric
spaces with distances $d$ and $d^{\,\prime},$ respectively. A family
$\frak{F}$ of continuous mappings from $X$ into $X^{\,\prime}$ is
said to be a {\bf normal} if every sequence of mappings $f_m$ in
$\frak{F}$ has a subsequence $f_{m_k}$ converging to a continuous
mapping $f:X\to X'$ uniformly on each compact set $C\subset X$.
Normality is closely related to the following notion. A family
$\frak{F}$ of mappings $f:X\rightarrow {X}^{\,\prime}$  is said to
be {\bf equicontinuous at a point} $x_0 \in X$ if for every
$\varepsilon > 0$ there is $\delta > 0$ such that ${d}^{\,\prime}
\left(f(x),f(x_0)\right)<\varepsilon$ for all $f \in \frak{F}$ and
$x \in X$ with $d(x,x_0)<\delta$. The family $\frak{F}$ is called
{\bf equicontinuous} if $\frak{F}$ is equicontinuous at every point
$x_0 \in X.$ The following version of the Arzela -- Ascoli theorem
will be useful later on, see e.g. Section 20.4 in \cite{Va}.

\bigskip
\begin{proposition}\label{pr3**!}{\it\, Let $(X, d)$ be a separable metric space and
let $\left(X^{\,\prime}, d^{\,\prime}\right)$ be a compact metric
space. Then a family $\frak{F}$ of mappings $f:X\rightarrow
X^{\,\prime}$ is normal if and only if $\frak{F}$ is
equicontinuous.}
\end{proposition}

\bigskip

In particular, Proposition \ref{pr3**!} holds in the case when
$X={{\Bbb R}}^n$ with the usual distance and $X'$ is the extended
space $\overline{{{\Bbb R}}^n}={{\Bbb R}}^n\bigcup\{\infty\}$
(compact) with the spherical metric. Recall that the {\bf spherical
(chordal) metric} $h(x,y)$ in $\overline{{{\Bbb R}}^n}$ is equal to
$|\pi(x)-\pi(y)|$ where $\pi$ is the stereographic projection of
$\overline{{{\Bbb R}}^n}$ on the sphere
$S^n(\frac{1}{2}e_{n+1},\frac{1}{2})$ in ${{\Bbb R}}^{n+1},$ i.e.,
in the explicit form,
$$h(x,\infty)=\frac{1}{\sqrt{1+{|x|}^2}}, \ \
h(x,y)=\frac{|x-y|}{\sqrt{1+{|x|}^2} \sqrt{1+{|y|}^2}}\,, \ \  x\ne
\infty\ne y\,.$$
The {\bf spherical diameter of a set} $E$ in $\overline{{\Bbb R}^n}$
is the quantity $h(E)=\sup\limits_{x_1, x_2\in E} h(x_1, x_2).$

\noindent Let $\frak{R}_{Q, \Delta}(D)$ be the class of all ring
$Q-$homeomorphisms $f$ in a domain $D\subseteq{\Bbb R}^n,$ $n\ge 2,$
such that $h\left(\overline{{\Bbb R}^n} \backslash f(D)\right) \ge
\Delta
>0.$ The following distortion estimate
under $Q$--ho\-me\-o\-mor\-phisms can be found in \cite{RS}, see
also Theorem 7.3 in $\cite{MRSY}.$

\bigskip
\begin{proposition}\label{pr2*!!}{\it\, Let $\Delta>0,$ $Q:D \rightarrow [0,\infty]$ be a
measurable function. Then
%for
\begin{equation}\label{eq2*!!}
h\left(f(x), f(x_0)\right)\le\frac{\alpha_n}{\Delta}
\exp\left\{-\int\limits_{|x-x_0|}^{\varepsilon(x_0)}
\frac{dr}{rq_{x_0}^{\frac{1}{n-1}}(r)}\right\}
\end{equation}
for every $f\in \frak{R}_{Q, \Delta}(D)$ and $x\in B(x_0,
\varepsilon(x_0)),$ $\varepsilon(x_0)< {\rm dist\,} \left(x_0,
\partial D\right),$ where $\alpha_n>0$ depends only on $n$ and
$q_{x_0}(r)$ is the mean value of the function $Q$ over the sphere
$|z-x_0|=r$.}
\end{proposition}
\medskip

For every non-decreasing function $\Phi:[0,\infty ]\to [0,\infty ]
,$ the {\bf inverse function} $\Phi^{-1}:[0,\infty ]\to [0,\infty ]$
can be well defined by setting
\begin{equation}\label{eq5.5CC} \Phi^{-1}(\tau)\ =\
\inf\limits_{\Phi(t)\ge \tau}\ t\ .
\end{equation} As usual, here $\inf$ is equal to $\infty$ if the set of
$t\in[0,\infty ]$ such that $\Phi(t)\ge \tau$ is empty. Note that
the function $\Phi^{-1}$ is non-decreasing, too.\medskip

\begin{remark}\label{rmk3.333} Immediately by the
definition  it is evident  that
\begin{equation}\label{eq5.5CCC} \Phi^{-1}(\Phi(t))\ \le\ t\ \ \ \ \
\ \ \ \forall\ t\in[ 0,\infty ]
\end{equation} with the equality in (\ref{eq5.5CCC}) except
intervals of constancy of the function $\Phi(t)$.
\end{remark}

\medskip
Since the mapping $t\mapsto t^p$ for every positive $p$ is a
sense--preserving homeomorphism $[0, \infty]$ onto $[0, \infty]$ we
may rewrite Theorem 2.1 from \cite{RSY$_1$} in the following form
which is more convenient for further applications. Here, in
(\ref{eq333Y}) and (\ref{eq333F}), we complete the definition of
integrals by $\infty$ if $\Phi_p(t)=\infty ,$ correspondingly,
$H_p(t)=\infty ,$ for all $t\ge T\in[0,\infty) .$ The integral in
(\ref{eq333F}) is understood as the Lebesgue--Stieltjes integral and
the integrals  in (\ref{eq333Y}) and (\ref{eq333B})--(\ref{eq333A})
as the ordinary Lebesgue in\-te\-grals.

\medskip
\begin{proposition} \label{pr4.1aB}{\it\, Let $\Phi:[0,\infty ]\to [0,\infty ]$ be a
non-decreasing function. Set \begin{equation}\label{eq333E} H_p(t)\
=\ \log \Phi_p(t)\ , \qquad \Phi_p(t)=\Phi\left(t^p\right)\,,\quad
p\in (0,\infty)\,.\end{equation}

Then the equality \begin{equation}\label{eq333Y}
\int\limits_{\delta}^{\infty} H^{\,\prime}_p(t)\ \frac{dt}{t}\ =\
\infty
\end{equation} implies the equality \begin{equation}\label{eq333F}
\int\limits_{\delta}^{\infty} \frac{dH_p(t)}{t}\ =\ \infty
\end{equation} and (\ref{eq333F}) is equivalent to
\begin{equation}\label{eq333B} \int\limits_{\delta}^{\infty}H_p(t)\
\frac{dt}{t^2}\ =\ \infty
\end{equation} for some $\delta>0,$ and (\ref{eq333B}) is equivalent to
every of the equalities: \begin{equation}\label{eq333C}
\int\limits_{0}^{\Delta}H_p\left(\frac{1}{t}\right)\ {dt}\ =\ \infty
\end{equation} for some $\Delta>0,$ \begin{equation}\label{eq333D}
\int\limits_{\delta_*}^{\infty} \frac{d\eta}{H_p^{-1}(\eta)}\ =\
\infty
\end{equation} for some $\delta_*>H(+0),$ \begin{equation}\label{eq333A}
\int\limits_{\delta_*}^{\infty}\ \frac{d\tau}{\tau \Phi_p^{-1}(\tau
)}\ =\ \infty \end{equation} for some $\delta_*>\Phi(+0).$
\medskip

Moreover, (\ref{eq333Y}) is equivalent  to (\ref{eq333F}) and hence
(\ref{eq333Y})--(\ref{eq333A})
 are equivalent each to other  if $\Phi$ is in addition absolutely continuous.
In particular, all the conditions (\ref{eq333Y})--(\ref{eq333A}) are
equivalent if $\Phi$ is convex and non--decreasing.}
\end{proposition}

\medskip
It is easy to see that conditions (\ref{eq333Y})--(\ref{eq333A})
become weaker as $p$ increases, see e.g. (\ref{eq333B}). It is
necessary to give one more explanation. From the right hand sides in
the conditions (\ref{eq333Y})--(\ref{eq333A}) we have in mind
$+\infty$. If $\Phi_p(t)=0$ for $t\in[0,t_*]$, then $H_p(t)=-\infty$
for $t\in[0,t_*]$ and we complete the definition $H_p'(t)=0$ for
$t\in[0,t_*]$. Note, the conditions (\ref{eq333F}) and
(\ref{eq333B}) exclude that $t_*$ belongs to the interval of
integrability because in the contrary case the left hand sides in
(\ref{eq333F}) and (\ref{eq333B}) are either equal to $-\infty$ or
indeterminate. Hence we may assume in (\ref{eq333Y})--(\ref{eq333C})
that $\delta>t_0$, correspondingly, $\Delta<1/t_0$ where $t_0\colon
=\sup\limits_{\Phi_p(t)=0}t$, $t_0=0$ if $\Phi_p(0)>0$.

\medskip
\section{The main lemma and its corollaries}
\setcounter{equation}{0}

Recall that a function  $\Phi :[0,\infty ]\to [0,\infty ]$ is called
{\bf convex} if
$$
\Phi (\lambda t_1 + (1-\lambda) t_2)\ \le\ \lambda\ \Phi (t_1)\ +\
(1-\lambda)\ \Phi (t_2)$$ for all $t_1$ and $t_2\in[0,\infty ] $ and
$\lambda\in [0,1]$.\medskip

In what follows, ${\Bbb R}^n(\varepsilon),$ $\varepsilon\in (0,1)$
denotes the ring in the space ${\Bbb R}^n$, $n\ge 2$,
\begin{equation}\label{eq5.5Cf} {\Bbb R}^n(\varepsilon)\ =\ \{\ x\in{\Bbb R}^n:\
\varepsilon<|x|\ <\ 1\ \}\ .\end{equation}

\noindent The following statement is a generalization and
strengthening of Lemma 3.1 from \cite{RSY$_1$}.

\begin{lemma} \label{lem5.5C} {\it\, Let $Q:{\Bbb B}^n\to [0,\infty ]$ be a measurable
function and let $\Phi:[0,\infty ]\to (0,\infty ]$ be a
non-decreasing convex function. Suppose that the mean value
$M(\varepsilon)$ of the function $\Phi\circ Q$ over the ring ${\Bbb
R}^n(\varepsilon),$ $\varepsilon\in (0, 1),$ is finite. Then
\begin{equation}\label{eq3.222} \int\limits_{\varepsilon}^{1}\
\frac{dr}{rq^{\frac{1}{p}}(r)}\ \ge\ \frac{1}{n}\
\int\limits_{eM(\varepsilon)}^{\frac{M(\varepsilon)}{\varepsilon^n}}\
\frac{d\tau}{\tau \left[\Phi^{-1}(\tau
)\right]^{\frac{1}{p}}}\qquad\qquad\forall\quad p\in (0, \infty)
\end{equation} where $q(r)$ is the average of the function $Q(x)$
over the sphere $|x|=r.$ }\end{lemma}

\bigskip

\begin{remark}\label{rmk3.333A} Note that (\ref{eq3.222}) is
equivalent  for each $p\in (0, \infty)$ to the inequality
\begin{equation}\label{eq3.1!}
\int\limits_\varepsilon^1\frac{dr}{rq^{\frac{1}{p}}(r)}\ \ge\
\frac{1}{n}\int\limits_{eM(\varepsilon)}^{\frac{M(\varepsilon)}{\varepsilon^n}}\frac{d\tau}{\tau\Phi_p^{\,-1}(\tau)}\
,\qquad \Phi_p(t)\ \colon =\ \Phi(t^p)\ .
\end{equation} Note also that $M(\varepsilon)$ converges as $\varepsilon\to
0$ to the average of $\Phi\circ Q$ over the unit ball ${\Bbb B}^n$.
\end{remark}

{\it Proof.} Denote $t_*=\sup\limits_{\Phi_p(t)=\tau_0} t,$
$\tau_0=\Phi(0).$ Setting $H_p(t) = \log\ \Phi_p(t),$ we see that
$H_p^{-1}(\eta)\ =\ \Phi_p^{-1}(e^{\eta})\ ,\ \ \ \Phi_p^{-1}(\tau)\
=\ H_p^{-1}(\log\ \tau)\ .$ Thus, we obtain that
$$q^{\frac{1}{p}}(r) = H_p^{-1}\left(\log \frac{h(r)}{r^n}\right) =
H_p^{-1}\left(n\log \frac{1}{r} + \log\ h(r)\right)\ \ \ \ \ \ \ \
\forall\ r\ \in R_*$$ where $h(r) \colon=\ r^n\Phi\left(q(r)\right)\
=\ r^n\Phi_p\left(q^{\frac{1}{p}}(r)\right)$ and $R_* = \{\
r\in(\varepsilon, 1):\ q^{\frac{1}{p}}(r)\
>\ t_*\}$. Then also
\begin{equation}\label{eq5.555K}q^{\frac{1}{p}}(e^{-s})\ =\ H_p^{-1}\left(ns\ +\ \log\
h(e^{-s})\right)\ \ \ \ \ \ \ \ \forall\ s\ \in S_*\end{equation}
where $ S_*\ =\ \{ s\in(0, \log\frac{1}{\varepsilon}):\
q^{\frac{1}{p}}\left(e^{-s}\right)\
>\ t_*\}$.

Now, by the Jensen inequality and convexity of $\Phi$ we have that
$$ \int\limits_0^{\log \frac{1}{\varepsilon}} h(e^{-s})\ ds\ =\
\int\limits_{\varepsilon}^{1} h(r)\ \frac{dr}{r}\ =\
\int\limits_{\varepsilon}^{1} \Phi(q(r))\ r^{n-1}{dr}$$
$$\le\ \int\limits_{\varepsilon}^{1}\left(\dashint_{S(r)}
\Phi(Q(x))\ d{\cal {A}}\right)\ r^{n-1}{dr}\ \le
\frac{\Omega_{n}}{\omega_{n-1}}\cdot M(\varepsilon)=\frac{1}{n}\cdot
M(\varepsilon)$$ where we use the mean value of the function
$\Phi\circ Q$ over the sphere $S(r)=\{ x\in{\Bbb R}^n: |x|=r\}$ with
respect to the area measure. As usual, here $\Omega_{n}$ and
$\omega_{n-1}$ is the volume of the unit ball and the area of the
unit sphere in ${\Bbb R}^n,$  correspondingly. Then arguing by
contradiction it is easy to see that
\begin{equation}\label{eq5.555N} |T|\ =\ \int\limits_{T}ds\ \le\
\frac{1}{n}\end{equation} where $T\ =\ \{\ s\in (0,
\log\frac{1}{\varepsilon}):\ \ \ h(e^{-s})\ > M(\varepsilon)\}.$
Next, let us show that
\begin{equation}\label{eq5.555O}q^{\frac{1}{p}}\left(e^{-s}\right)
\ \le\ H_p^{-1}\left(ns\ +\ \log\ M(\varepsilon)\right)\ \ \ \ \ \ \
\ \ \ \forall\ s\in\left(0, \log\frac{1}{\varepsilon}\right)
\setminus  T_*\end{equation} where $T_*\ =\ T\cap S_*$. Note that
$\left(0,\log\frac{1}{\varepsilon}\right)\setminus T_*  =
\left[\left(0,\log\frac{1}{\varepsilon}\right)\setminus S_*\right]
\cup \left[\left(0, \log\frac{1}{\varepsilon}\right)\setminus
T\right] = \left[\left(0, \log\frac{1}{\varepsilon}\right)\setminus
S_*\right] \cup \left[S_*\setminus T\right]$. The inequality
(\ref{eq5.555O}) holds for $s\in S_*\setminus T$ by (\ref{eq5.555K})
because $H_p^{-1}$ is a non-decreasing function. Note also that
$e^{ns}M(\varepsilon)> \Phi(0) = \tau_0$ for all $s\in\left(0,\
\log{1/\varepsilon}\right)$ and then $t_* <
\Phi_p^{-1}\left(e^{ns}M(\varepsilon) \right) = H_p^{-1}\left(ns\ +\
\log\ M(\varepsilon)\right)$ for all $s\in\left(0,\ \log
{1/\varepsilon}\right)\ .$ Consequently, (\ref{eq5.555O}) holds for
$s\in(0, \log\frac{1}{\varepsilon})\setminus S_*$, too. \medskip

Since $H_p^{-1}$ is non--decreasing, we have by (\ref{eq5.555N}) and
(\ref{eq5.555O}) that \begin{equation}\label{eq5.555P}
\int\limits_{\varepsilon}^{1}\ \frac{dr}{rq^{\frac {1}{p}}(r)}\ =\
\int\limits_{0}^{\log \frac{1}{\varepsilon}}\ \frac{ds}{q^\frac
{1}{p}(e^{-s})}\ \ge\
\int\limits_{\left(0,\log\frac{1}{\varepsilon}\right)\setminus T_*}\
\frac{ds}{H_p^{-1}(ns + \Delta)}\ \ge\ \end{equation}
$$
\ge\ \int\limits_{|T_*|}^{\log \frac{1}{\varepsilon}}\
\frac{ds}{H_p^{-1}(ns + \Delta)}\ \ge\
\int\limits_{\frac{1}{n}}^{\log \frac{1}{\varepsilon}}\
\frac{ds}{H_p^{-1}(ns + \Delta)}\ =\
\frac{1}{n}\int\limits_{1+{\Delta}}^{n\log
\frac{1}{\varepsilon}+\Delta}\ \frac{d\eta}{H_p^{-1}(\eta)}
$$
where $ \Delta=\log {M(\varepsilon)}$. Note that $1+\Delta = \log\
e{M(\varepsilon)}.$ Thus,
\begin{equation}\label{eq5.555S} \int\limits_{\varepsilon}^{1}\
\frac{dr}{rq^{\frac{1}{p}}(r)}\ \ge\ \frac{1}{n}\int\limits_{\log
e{M(\varepsilon)} }^{\log\frac{M(\varepsilon)}{\varepsilon^n}}\
\frac{d\eta}{H_p^{-1}(\eta)}
\end{equation} and, after the replacement $\eta = \log\ \tau$, we
obtain (\ref{eq3.1!}) and hence (\ref{eq3.222}).

\medskip
\begin{corollary} \label{cor3.1}{\,\it
Let $\Phi:[0,\infty ]\rightarrow (0,\infty ]$ be a non-decreasing
convex function, $Q:{\Bbb B}^n\rightarrow [0,\infty ]$ a measurable
function, $Q_*(x)=1$ if $Q(x)<1$ and $Q_*(x)=Q(x)$ if $Q(x)\ge 1$.
Suppose that the mean $M_*(\varepsilon)$ of the function $\Phi\circ
Q$ over the ring ${\Bbb R}^n(\varepsilon),$ $\varepsilon\in (0, 1),$
is finite. Then
\begin{equation}\label{eq3.1}
\int\limits_{\varepsilon}^{1}\ \frac{dr}{rq^{\frac{\lambda}{p}}(r)}\
\ge\ \frac{1}{n}\
\int\limits_{eM_*(\varepsilon)}^{\frac{M_*(\varepsilon)}{\varepsilon^n}}\
\frac{d\tau}{\tau \left[\Phi^{-1}(\tau )\right]^{\frac{1}{p}}}\ \ \
\qquad \ \ \ \ \ \ \ \forall\ \lambda\ \in\ (0,1), \qquad p\in (0,
\infty)
\end{equation}
where $q(r)$ is the average of the function $Q(x)$ over the sphere
$|x|=r.$ }
\end{corollary}

\medskip
\medskip
Indeed, let $q_*(r)$ be the average of the function $Q_*(x)$ over
the sphere $|x|=r$. Then $q(r)\le q_*(r)$ and, moreover, $q_*(r)\ge
1$ for all $r\in (0,1)$. Thus, $q^{\frac{\lambda}{p}}(r)\le
q_*^\frac{\lambda}{p}(r)\le q_*^\frac{1}{p}(r)$  for all $\lambda\in
(0,1)$ and hence by Lemma \ref{lem5.5C} applied to $Q_*(x)$ we
obtain (\ref{eq3.1}).

\medskip
\begin{theorem} \label{th5.555}{\it\, Let $Q:{\Bbb B}^n\to [0,\infty ]$ be a measurable
function such that \begin{equation}\label{eq5.555}
\int\limits_{{\Bbb B}^n} \Phi (Q(x))\ dm(x)\  <\
\infty\end{equation} where $\Phi:[0,\infty ]\to [0,\infty ]$ is a
non-decreasing convex function such that
\begin{equation}\label{eq3.333a} \int\limits_{\delta_0}^{\infty}\ \frac{d\tau}{\tau
\left[\Phi^{-1}(\tau )\right]^{\frac{1}{p}}}\ =\ \infty\,,\qquad
p\in (0, \infty)\,,
\end{equation} for some $\delta_0\
>\ \tau_0\ \colon =\ \Phi(0).$ Then \begin{equation}\label{eq3.333A}
\int\limits_{0}^{1}\ \frac{dr}{rq^{\frac{1}{p}}(r)}\ =\ \infty
\end{equation} where $q(r)$ is the average of the function $Q(x)$
over the sphere $|x|=r$.}
\end{theorem}

\begin{remark}\label{rmk4.7www} Since $\left[\Phi^{\,-1}(\tau)\right]^{\frac{1}{p}}=
\Phi_p^{\,-1}(\tau)$ where $\Phi_p(t)=\Phi(t^p),$  (\ref{eq3.333a})
implies that
\begin{equation}\label{eq3.a333} \int\limits_{\delta}^{\infty}\ \frac{d\tau}{\tau
\Phi^{-1}_p(\tau )}\ =\ \infty\ \ \ \ \ \ \ \ \ \ \forall\ \delta\
\in\ [0,\infty)   \end{equation} but (\ref{eq3.a333}) for some
$\delta\in[0,\infty)$, generally speaking, does not imply
(\ref{eq3.333a}). Indeed, for $\delta\in [0,\delta_0),$
(\ref{eq3.333a}) evidently implies (\ref{eq3.a333}) and, for
$\delta\in(\delta_0,\infty)$, we have that
\begin{equation}\label{eq3.e333} 0\ \le\ \int\limits_{\delta_0}^{\delta}\
\frac{d\tau}{\tau \Phi_p^{-1}(\tau )}\ \le\
\frac{1}{\Phi_p^{-1}(\delta_0)}\ \log\ \frac{\delta}{\delta_0}\ <\
\infty
\end{equation} because $\Phi_p^{-1}$ is non-decreasing and
$\Phi_p^{-1}(\delta_0)>0$. Moreover, by the de\-fi\-ni\-tion of the
inverse function $\Phi_p^{-1}(\tau)\equiv 0$ for all $\tau \in
[0,\tau_0],$ $\tau_0=\Phi_p(0)$, and hence (\ref{eq3.a333}) for
$\delta\in[0,\tau_0),$ generally speaking, does not imply
(\ref{eq3.333a}). If $\tau_0 > 0$, then
\begin{equation}\label{eq3.c333} \int\limits_{\delta}^{\tau_0}\
\frac{d\tau}{\tau \Phi_p^{-1}(\tau )}\ =\ \infty\ \ \ \ \ \ \ \ \ \
\forall\ \delta\ \in\ [0,\tau_0)  \end{equation} However,
(\ref{eq3.c333}) gives no information on the function $Q(x)$ itself
and, consequently, (\ref{eq3.a333}) for $\delta < \Phi(0)$ cannot
imply (\ref{eq3.333A}) at all. \end{remark}

\medskip
In view of (\ref{eq3.a333}), Theorem \ref{th5.555} follows
immediately from Lemma \ref{lem5.5C}.

\medskip
\begin{corollary} \label{cor555}{\it\, If $\Phi:[0,\infty ]\to [0,\infty ]$ is a
non-decreasing convex func\-tion and $Q$ satisfies the condition
(\ref{eq5.555}), then each of the conditions
(\ref{eq333Y})--(\ref{eq333A}) for $p\in (0, \infty)$ implies
(\ref{eq3.333A}). Moreover, if in addition $\Phi(1)<\infty$ or
$q(r)\ge 1$ on a subset of $(0,1)$ of a positive measure, then each
of the conditions (\ref{eq333Y})--(\ref{eq333A}) for $p\in (0,
\infty)$ implies
\begin{equation}\label{eq3.3} \int\limits_{0}^{1}\
\frac{dr}{rq^{\frac{\lambda}{p}}(r)}\ =\ \infty\ \ \ \ \ \ \ \ \
\forall\ \lambda\ \in\ (0,1)
\end{equation}
and also
\begin{equation}\label{eq3.3AB} \int\limits_{0}^{1}\
\frac{dr}{r^{\alpha}q^{\frac{\beta}{p}}(r)}\ =\ \infty\ \ \ \ \ \ \
\ \ \forall\ \alpha\ge 1 ,\ \beta\ \in\ (0,\alpha]
\end{equation}}
\end{corollary}

\section{Sufficient conditions for equicontinuity}
\setcounter{equation}{0}

\medskip
Let $D$ be a fixed domain in the extended space $\overline{{\Bbb
R}^n}={\Bbb R}^n\cup\{\infty\},$ $n\ge 2.$ Given a function
$\Phi:[0, \infty]\rightarrow [0, \infty],$ $M>0,$ $\Delta>0$,
$\frak{R}^{\Phi}_{M,\Delta}$ denotes the collection of all ring
$Q(x)$--ho\-me\-o\-mor\-phisms in $D$ such that
$h\left(\overline{{\Bbb R}^n}\setminus f(D)\right)\ge \Delta$ and
\begin{equation}\label{eq2!!}
\int\limits_D\Phi\left(Q(x)\right)\frac{dm(x)}{\left(1+|x|^2\right)^n}\
\le\ M\,.
\end{equation}

\medskip
\begin{theorem}\label{th1!}{\it\, Let
$\Phi:[0, \infty]\rightarrow [0, \infty]$ be non-decreasing convex
function. If
\begin{equation}\label{eq3!}
\int\limits_{\delta_0}^{\infty}
\frac{d\tau}{\tau\left[\Phi^{-1}(\tau)\right]^{\frac{1}{n-1}}}\ =\
\infty
\end{equation}
for some $\delta_0>\tau_0:=\Phi(0),$ then the class
$\frak{R}^{\Phi}_{M,\Delta}$ is equicontinuous and, consequently,
forms a normal family of mappings for every $M\in(0, \infty)$ and
$\Delta\in(0, 1).$ }
\end{theorem}

\medskip
\begin{remark}\label{rem1}
Note that the condition
\begin{equation}\label{eq3!!}
\int\limits_D \Phi\left(Q(x)\right)dm(x)\le M
\end{equation}
implies (\ref{eq2!!}). Thus, the condition (\ref{eq2!!}) is more
general than (\ref{eq3!!}) and ring $Q$--homeomorphisms satisfying
(\ref{eq3!!}) form a subclass of $\frak{R}^{\Phi}_{M,\Delta}.$
Conversely, if the domain $D$ is bounded, then (\ref{eq2!!}) implies
the condition
\begin{equation}\label{eq4!}
\int\limits_D \Phi\left(Q(x)\right)dm(x)\le M_*
\end{equation}
where $M_*=M\cdot\left(1+\delta_*^2\right),$
$\delta_*=\sup\limits_{x\in D}|x|.$
\end{remark}

\medskip
\medskip
{\it Proof.} With no loss of generality we may assume that
$\Phi(0)>0.$ By Proposition \ref{pr3**!} it is sufficient to show
that mappings in $\frak{R}^{\Phi}_{M,\Delta}$ are equicontinuous at
every point $x_0\in D.$ If $x_0\ne \infty,$ then by the Proposition
\ref{pr2*!!}
\begin{equation}\label{eq6!}
h\left(f(x), f(x_0)\right)\ \le\ \frac{\alpha_n}{\Delta}
\exp\left\{-\int\limits_{|x-x_0|}^{\rho}
\frac{dr}{rq_{x_0}^{\frac{1}{n-1}}(r)}\right\}
\end{equation}
for every fixed $x\in B(x_0, \rho)$ and every positive
$\rho=\rho(x_0)< {\rm dist\,}(x_0, \partial D)$ where $q_{x_0}(r)$
is the mean value of $Q(x)$ over the sphere $|z-x_0|=r$ and
$\alpha_n$ depends only on $n.$ After the replacement
$t={r}/{\rho},$ we have that the integral from the right hand side
in (\ref{eq6!}) is estimated by Lemma \ref{lem5.5C} in the following
way
$$\int\limits_{|x-x_0|}^{\rho}\frac{dr}{rq_{x_0}^{\frac{1}{n-1}}(r)}
=\int\limits_{\varepsilon}^1\frac{dt}{tq^{\frac{1}{n-1}}(t)}\ge
\frac{1}{n}\int\limits_{eM(\varepsilon)}^{\frac{M(\varepsilon)}{\varepsilon^n}}\frac{d\tau}
{\tau\left[\Phi^{-1}(\tau)\right]^{\frac{1}{n-1}}}$$
where $\varepsilon=|x-x_0|/\rho,$ $q(t)=q_{x_0}(\rho t)$ and
$$M(\varepsilon)=\dashint_R \Phi\left(Q(z)\right)dm(z)=
\frac{1}{\Omega_n\rho^n\left(1-\varepsilon^n\right)}\int\limits_R
\Phi\left(Q(z)\right)dm(z)$$
where $R=\left\{z\in {\Bbb R}^n: |x-x_0|<|z-x_0|<\rho\right\}$ is a
ring centered at $x_0$ and $\Omega_n$ is the volume of the unit ball
${\Bbb B}^n$ in ${\Bbb R}^n.$ Note that
$$M(\varepsilon)\le \frac{\beta_n(x_0)}{\Omega_n
(1-\varepsilon^n)}\int\limits_R
\Phi(Q(z))\frac{dm(z)}{\left(1+|z|^2\right)^n}$$
where $\beta_n(x_0)=\left(1+(\rho(x_0)+|x_0|)^2\right)^n
/\rho^n(x_0)$ because $|z|\le |z-x_0|+|x_0|\le \rho(x_0)+|x_0|.$
Thus,
$$\Phi(0)\le M(\varepsilon)\le \frac{2\beta_n(x_0)}{\Omega_n}M$$
if $\varepsilon\le 1/\sqrt[n]{2}$ and, in particular, if
$\varepsilon\le 1/2.$ Consequently,
\begin{equation}\label{eq7!}
h\left(f(x), f(x_0)\right)\quad\le\quad\frac{\alpha_n}{\Delta}
\exp\left\{-\frac{1}{n}\int\limits_{\lambda_n\beta_n(x_0)M}^{\frac{\Phi(0)\rho^n(x_0)}{|x-x_0|^n}}
\frac{d\tau}{\tau \left[\Phi^{-1}(r)\right]^{\frac{1}{n-1}}}\right\}
\end{equation}
for all $x$ such that $|x-x_0|<\rho(x_0)/2$ where
$\lambda_n=2e/\Omega_n$ depends only on $n.$ Thus, $f\in
\frak{R}_{M, \Delta}^{\Phi}$ are equicontinuous at the point $x_0.$
The case $x_0=\infty$ is reduced to $x_0=0$ by the inversion with
respect to the unit sphere $|x|=1$.

\medskip
\begin{corollary}\label{cor1!}{\,\it
Each of the conditions (\ref{eq333Y})--(\ref{eq333A}) for $p\in (0,
n-1] $ implies equicontinuity and normality of the classes
$\frak{R}^{\Phi}_{M,\Delta}$ for all $M\in (0, \infty)$ and
$\Delta\in (0, 1).$ }
\end{corollary}

\medskip
Given a function $\Phi:[0, \infty]\rightarrow [0, \infty],$ $M>0$
and $\Delta>0,$ $S^{\Phi}_{M, \Delta}$ denotes the class of all
homeomorphisms $f$ of $D$ in the Sobolev class $W_{loc}^{1,n}$ with
a locally integrable $K_I(x,f)$ such that $h\left(\overline{{\Bbb
R}^n}\setminus f(D)\right)\ge\Delta$ and (\ref{eq2!!}) holds for
$Q(x)=K_I(x,f).$ Note that if $\Phi$ is non-decreasing, convex and
non--constant on $[0,\infty)$, then (\ref{eq2!!}) itself implies
that $K_I(x,f)\in L_{loc}^1.$ Note also that $S^{\Phi}_{M,
\Delta}\subset \frak{R}^{\Phi}_{M, \Delta},$ see e.g. Theorem 4.1 in
\cite{MRSY}. Thus, we have the following consequence.
\medskip

\begin{corollary}\label{cor2!}{\,\it
Each of the conditions (\ref{eq333Y})--(\ref{eq333A}) for $p\in (0,
n-1]$ implies equicontinuity and normality of the class
$S^{\Phi}_{M,\Delta}$ for all $M\in (0, \infty)$ and $\Delta\in (0,
1).$ }
\end{corollary}

\medskip
\begin{remark}\label{rem2}
The given conditions (\ref{eq333Y})--(\ref{eq333A}) for $p=n-1$ are
weakest that lead to equicontinuity (normality) of the classes
$S_{M, \Delta}^{\Phi}$ and $\frak{R}_{M, \Delta}^{\Phi},$ see
Theorem \ref{th3} further. The most interesting of them is
(\ref{eq333B}) that can be rewritten in the following form:
\begin{equation}\label{eq5!}
\int\limits_{\delta}^{\infty}\log \Phi(t)\ \
\frac{dt}{t^{n^{\,\prime}}}=\infty
\end{equation}
where $\frac{1}{n^{\,\prime}}+\frac{1}{n}=1,$ i.e. $n^{\,\prime}=2$
for $n=2,$ $n^{\,\prime}$ is strictly increasing in $n$ and
$n^{\prime}=n/(n-1)\rightarrow 1$ as $n\rightarrow \infty.$ Note
also that the condition (\ref{eq3!}), as well as (\ref{eq3})
further, can be rewritten in the form \begin{equation}\label{eq33}
\int\limits_{\delta}^{\infty}\frac{d\tau}{\tau
\Phi_{n-1}^{\,-1}(\tau)}\ =\ \infty\ ,\ \ \ \ \ \ \ \ \
\Phi_{n-1}(t)\ \colon =\ \Phi(t^{n-1})\ .
\end{equation}
\end{remark}

\section{Necessary conditions for equicontinuity}
\setcounter{equation}{0}

\begin{theorem}\label{th3}{\it\, If the
classes $S^{\Phi}_{M,\Delta}\subset \frak{R}^{\Phi}_{M,\Delta}$ are
equicontinuous (normal) for a non--decreasing convex function
$\Phi:[0, \infty]\rightarrow [0, \infty],$ all $M\in (0,\infty)$ and
$\Delta\in (0,1).$ Then
\begin{equation}\label{eq3}
\int\limits_{\delta_*}^{\infty}\frac{d\tau}{\tau
\left[\Phi^{\,-1}(\tau)\right]^{\frac{1}{n-1}}}\ =\ \infty
\end{equation}} for all $\delta_*\in (\tau_0,
\infty)$ where $\tau_0\ \colon=\ \Phi(0).$
\end{theorem}\medskip

It is evident that the function $\Phi(t)$ in Theorem \ref{th3}
cannot be constant because in the contrary case we would have no
real restrictions for $K_{I}$ except $\Phi(t)\equiv\infty$ when the
classes $S^{\Phi}_{M,\Delta}$ are empty. Moreover, by the known
criterion of convexity, see e.g. Proposition 5 in I.4.3 of
\cite{Bou}, the slope $[\Phi(t)-\Phi(0)]/t$ is nondecreasing. Hence
the proof of Theorem \ref{th3} follows from the next
statement.\medskip

\begin{lemma}\label{th3!}{\it\, Let a function $\Phi : [0,\infty]\to[0,\infty]$
be non-decreasing and %
\begin{equation}\label{eq4!!}
\Phi(t)\ \ge\ C\cdot t^{\frac{1}{n-1}}\qquad\forall\ t\in [T,
\infty]
\end{equation}
for some $C>0$ and $T\in (0, \infty).$ If the classes
$S_{M,\Delta}^{\Phi}\subset \frak{R}_{M,\Delta}^{\Phi}$ are
equicontinuous (normal) for all $M\in (0,\infty)$ and $\Delta\in
(0,1)$, then (\ref{eq3}) holds for all $\delta_*\in (\tau_0,
\infty)$ where $\tau_0\ \colon=\ \Phi(+0).$}
\end{lemma}

\medskip
\begin{remark}\label{rem3} As well--known, the critical exponent
$n-1$ takes a key part in many problems of space mappings. The
condition (\ref{eq4!!}) can be rewritten in the form
\begin{equation}\label{eq4!!!}
\Phi_{n-1}(t)\ \ge\ C\cdot t\qquad\forall\ t\in [T, \infty]
\end{equation}
where $\Phi_{n-1}(t)=\Phi(t^{n-1})$ and $C>0,$ $T\in (0, \infty)$
that once more accentuates the significance of the function
$\Phi_{n-1}$ in the question. In fact, it suffices also to require
the weaker condition of convexity of $\Phi_{n-1}$ instead of $\Phi$
in Theorem \ref{th3}.
\end{remark}
\medskip

{\it Proof of Lemma \ref{th3!}}. Let us assume that (\ref{eq3}) is
not true, i.e.
\begin{equation}\label{eq3!!!}
\int\limits_{\delta_0}^{\infty}\frac{d\tau}{\tau
\Phi_{n-1}^{\,-1}(\tau)}\ <\ \infty\end{equation}
for some $\delta_0\in (\tau_0, \infty)$ where $\ \Phi_{n-1}(t)\
\colon =\ \Phi(t^{n-1})\ .$ Then also
\begin{equation}\label{eq7!!}
\int\limits_{\delta}^{\infty}\frac{d\tau}{\tau
\Phi_{n-1}^{\,-1}(\tau)}\quad<\quad\infty\qquad \forall\quad
\delta\in (\tau_0, \infty)\end{equation}
because $\Phi^{-1}(\tau)>0$ for all $\tau>\tau_0$ and
$\Phi^{-1}(\tau)$ is non--decreasing. Note that by (\ref{eq4!!})
\begin{equation}\label{eq3.VVV!} \Phi_{n-1}(t)\ \ge\ {C}\cdot{t}\ \ \ \ \
\ \ \ \forall\ t\ \ge\ T \end{equation} under some $C>0$ and
$T\in(1,\infty)$. Furthermore, applying the linear transformation
$\alpha\Phi + \beta$ with $\alpha = 1/C$ and $\beta = T,$ see e.g.
(\ref{eq333B}), we may assume that
\begin{equation}\label{eqKKK3} \Phi_{n-1}(t)\ \ge\ t\ \ \ \ \ \ \ \ \ \ \ \ \forall\
t\in[0,\infty)\ .\end{equation}
Of course, we may also assume that $\Phi(t)=t$ for all $t\in[0,1)$
because the values of $\Phi$ in $[0,1)$ give no information on
$K_I(x,f)\ge 1$ in (\ref{eq2!!}). It is clear that (\ref{eq7!!})
implies $\Phi(t)<\infty$ for all $t<\infty,$ see the criterion
(\ref{eq333B}), cf. (\ref{eq333A}).

Now, note that the function $\Psi(t)\colon = t\Phi_{n-1}(t)$ is
strictly increasing, $\Psi(1)=\Phi(1)$ and $\Psi(t)\to\infty$ as
$t\to\infty$. Hence the functional equation
\begin{equation}\label{eqLLL3} \Psi(K(r))\ =\
\left(\frac{\gamma}{r}\right)^2\ \ \ \ \ \ \ \ \ \ \ \ \ \forall\ r\
\in\ (0,1]\ ,\end{equation} where $\gamma=\Phi^{1/2}(1)\ge 1$, is
well solvable with $K(1)=1$ and a strictly decreasing continuous
$K(r),$ $K(r)<\infty,$ $r\in (0, 1],$ and $K(r)\rightarrow \infty$
as $r\rightarrow 0.$ Taking the logarithm in (\ref{eqLLL3}), we have
that
$$
\log\ K(r)\ +\ \log\ \Phi_{n-1}(K(r))\ =\ 2\ \log\ \frac{\gamma}{r}
$$
and by (\ref{eqKKK3}) we obtain that
$$
\log\ K(r)\ \le\ \log\ \frac{\gamma}{r}\ ,
$$
i.e.,
\begin{equation}\label{eqMMM3} K(r)\ \le\ \frac{\gamma}{r}\ . \end{equation} Then by
(\ref{eqLLL3})
$$
\Phi_{n-1}(K(r))\ \ge\ \frac{\gamma}{r}
$$
and by (\ref{eq5.5CCC})
\begin{equation}\label{eq8!}
K(r)\ \ge\ \Phi_{n-1}^{-1}\left(\frac{\gamma}{r}\right)\ .
\end{equation}

It is sufficient to consider the case $D={\Bbb B}^n$. We define the
following mappings in the unit ball ${\Bbb B}^n$:
$$
f(x)=\frac{x}{|x|}\ R(|x|)\ ,\ \ \ \ \ \ \ \ \ f_m(x)=\frac{x}{|x|}\
R_m(|x|)\ ,\ \ \ m=1,2,\ldots
$$
where
$$
R(t)=\exp\{I(0)-I(t)\}\ ,\ \ \ \ \ \ R_m(t)=\exp\{I(0)-I_m(t)\}\ ,
$$
$$
I(t)=\int\limits_{t}^1\frac{dr}{rK(r)},
\qquad I_m(t)=\int\limits_{t}^1\frac{dr}{rK_m(r)}
$$
and
$$ K_m(r)\,=\,\left
\{\begin{array}{rr} K(r), & {\rm if } \ r\ge 1/m, \\
K\left(\frac{1}{m}\right), & {\rm if}\quad r\in (0, 1/m)\ .
\end{array}\right.
$$
By (\ref{eq8!})
$$
I(0)-I(t)=\ \int\limits_0^t\ \frac{dr}{rK(r)}\ \le\ \int\limits_0^t\
\frac{dr}{r\Phi_{n-1}^{-1}\left(\frac{\gamma}{r}\right)}\ =\
\int\limits_{\frac{\gamma}{t}}^{\infty}\
\frac{d\tau}{\tau\Phi_{n-1}^{-1}(\tau)} \qquad\forall\quad
t\in(0,1]$$
where $\gamma/t \ge \gamma \ge 1 > \Phi(0)=0.$ Hence by the
condition (\ref{eq7!!})
\begin{equation}\label{eqNNN3} I(0)-I(t)
\ \le\ I(0)\ =\ \int\limits_0^1\ \frac{dr}{rK(r)}\ <\ \infty\
\qquad\forall\quad t\in(0,1]
\end{equation}
Moreover, $f_m$ and $f\in C^1\left({\Bbb B}^n\setminus \{0\}\right)$
because $K_m(r)$ and $K(r)$ are continuous, and hence locally
quasiconformal in ${\Bbb B}^n\setminus \{0\}.$ Furthermore, $f_m$
are $K_m$--quasiconformal in ${\Bbb B}^n$ where
$K_m=K\left(1/m\right).$

Next, the tangent and radial distortions under the mapping $f$ on
the sphere $|x|=\rho,$ $\rho\in (0,1),$ are easy calculated
$$\delta_{\tau}(x)=
\frac{|f(x)|}{|x|}=\frac{\exp\left\{\int\limits_0^{\rho}\frac{dt}{K(t)}\right\}}
{\rho}\,,$$ $$\delta_r(x)=\frac{\partial |f(x)|}{\partial
|x|}=\frac{\exp\left\{\int\limits_0^{\rho}\frac{dt}{K(t)}\right\}}
{\rho K(\rho)}$$
and we see that $\delta_{\tau}(x)\le \delta_r(x)$ because $K(r)\ge
1.$
Consequently, by the spherical symmetry we have that
$$K_I(x,f)=\frac{\delta_{\tau}^{n-1}(x)\cdot\delta_r(x)}
{\delta^n_r(x)}=K^{n-1}(|x|)$$ at all points $x\in {\Bbb
B}^n\setminus\{0\},$ see e.g. Subsection I.4.1 in \cite{Re}. Note
that
\begin{equation}\label{eq5!!!}
f_m(x)\equiv f(x)\qquad \forall\ x:\ \frac{1}{m}<|x|< 1, \quad
m=1,2\ldots\,.
\end{equation}
Hence it is similarly calculated $K_I(x, f_m)=K_I(x,
f)=K^{n-1}(|x|)$ for $\frac{1}{m}<|x|<1$ and $K_I(x,f_m)=K(1/m)$ for
$0<|x|<\frac{1}{m}.$ Thus, $f_m$ are quasiconformal in ${\Bbb B}^n,$
hence $f_m\in W_{loc}^{1,n}$ and
by (\ref{eqLLL3})
$$
 \int\limits_{{\Bbb B}^n }\ \Phi\left(K_I(x,f_m)\right)\
 dm(x)\ \le\  \int\limits_{{\Bbb B}^n }\ \Phi_{n-1}\left(K(|x|)\right)\  dm(x)\ =$$
$$=\ \omega_{n-1}\int\limits_{0}^1\frac{\Psi\left(K(r)\right)}{rK(r)}\cdot
r^{n}dr\ \le\ \gamma^2\omega_{n-1}\int\limits_{0}^1\frac{dr}{rK(r)}\
\le\ M\ \colon =\ \gamma^2\omega_{n-1}I(0)\ <\ \infty\ .
$$
Note that $f_m$ map the unit ball ${\Bbb B}^n$ onto the ball
centered at the origin with the radius $e^{I(0)}<\infty$. Thus,
$f_m\in S^{\Phi}_{M,\Delta}$ with $M$ given above and some $\Delta >
0$.

On the other hand, it is easy to see that
\begin{equation}\label{eq2!}
\lim\limits_{x\rightarrow 0}\ |f(x)|\ =\ \lim\limits_{t\rightarrow
0}\ \rho(t)\ =\ e^{0}\ =\ 1\ ,
\end{equation}
i.e. $f$ maps the punctured ball ${\Bbb B}^n\setminus\{ 0\}$ onto
the ring $1<|y|< e^{I(0)}.$ Then by (\ref{eq5!!!}) and (\ref{eq2!})
we obtain that
$$|f_m(x)|=|f(x)|\ge 1\qquad\qquad\forall\quad x:|x|\ge 1/m,\quad m=1,2,\ldots\,,$$
i.e. the family $\{f_m\}_{m=1}^{\infty}$ is not equicontinuous at
the origin.

\medskip
The contradiction disproves the assumption (\ref{eq3!!!}).

\medskip
\begin{remark}\label{rem4}
Theorem \ref{th3} shows that the condition (\ref{eq3!}) in Theorem
\ref{th1!} is not only sufficient but also necessary for
equicontinuity (normality) of classes with the integral constraints
of the type either (\ref{eq2!!}) or (\ref{eq4!}) with a convex
non--decreasing $\Phi.$ In view of Proposition \ref{pr4.1aB}, the
same concerns to all the conditions (\ref{eq333Y})--(\ref{eq333A})
with $p=n-1.$
\end{remark}

\medskip
\begin{corollary}\label{cor3!}
{\it\, The equicontinuity (normality) of the classes
$S^{\Phi}_{M,\Delta}\subset \frak{R}^{\Phi}_{M,\Delta}$ for $M\in
(0, \infty)$, $\Delta\in (0,1)$ and non--decreasing convex $\Phi$
implies that
\begin{equation}\label{eq6!!}
\int\limits_{\delta}^{\infty}\log \Phi(t)\
\frac{dt}{t^{n^{\prime}}}\ =\ \infty
\end{equation}
for all $\delta>t_0$ where $t_0:=\sup\limits_{\Phi(t)=0}t,$ $t_0=0$
if $\Phi(0)>0,$
$\frac{1}{n^{\prime}}+\frac{1}{n}=1,$ i.e. $n^{\,\prime}=n/(n-1).$ }
\end{corollary}

\medskip
Recall that by Remark \ref{rem2} and Proposition \ref{pr4.1aB} the
condition (\ref{eq6!!}) is also sufficient for
equi\-con\-ti\-nui\-ty (normality) of the classes
$S^{\Phi}_{M,\Delta}$ and $\frak{R}^{\Phi}_{M,\Delta}$.

\bigskip
CONTACT INFORMATION
\bigskip

\noindent{{\it Vladimir Ryazanov and Evgenii Sevost'yanov}
\\{\sc Institute of Applied Mathematics and Mechanics \\
of National Academy of Sciences of Ukraine, \\
74  Roza Luksemburg Str., 83 114, Donetsk, Ukraine}, \\e--mail:
vlryazanov1@rambler.ru, brusin2006@rambler.ru

\end{document}